\documentclass[10pt]{amsart}
\usepackage{amsmath,amssymb}
\usepackage[letterpaper,margin=1in]{geometry}
\usepackage[hidelinks]{hyperref}
\newtheorem{theorem}{Theorem}
\newtheorem{proposition}[theorem]{Proposition}
\theoremstyle{remark}

\title[Connectivity of reduced bumpless pipe dreams]
  {A counterexample to connectivity under local flips of reduced bumpless pipe dreams}
\author{Gurleen Kaur Chahal \and Arpit Gupta}
\address[Gurleen Kaur Chahal]{School of Aeronautics, Neemrana, Rajasthan, India}
\address[Arpit Gupta]{Independent Researcher, Haridwar, Uttarakhand, India}
\email[Arpit Gupta]{arpitmathematics@gmail.com}
\thanks{Arpit Gupta is the corresponding author.}
\hypersetup{
  pdfauthor={Gurleen Kaur Chahal and Arpit Gupta},
  pdftitle={A counterexample to connectivity under local flips of reduced bumpless pipe dreams}
}

\begin{document}
\begin{abstract}
Anderson, Panova and Petrov conjectured that the set of reduced bumpless
pipe dreams of size \(n\) is connected under the \(2\times2\) flips that
preserve reducedness (Conjecture~4.2 of arXiv:2603.20104). We show that
the conjecture fails for \(n=9\): there are two disjoint sets of \(16\)
reduced bumpless pipe dreams, interchanged by transposition, each closed
under these flips and neither containing the identity (Rothe) pipe dream.
Exhaustive computation shows that the graph is connected for
\(n\leq 8\), so \(9\) is the smallest counterexample, and that the
remaining \(365{,}973{,}819\) pipe dreams of size \(9\) form a single
component. Padding the order-\(9\) example gives further components
avoiding the identity at orders \(10\), \(11\) and \(12\). The connectivity
theorem for flips augmented by droops and undroops, and the associated
MCMC sampler, are not affected.
\end{abstract}
\maketitle

\section{The conjecture and the counterexample}

A \emph{bumpless pipe dream} of size \(n\) is a tiling of the \(n\times n\)
grid by the six tiles of \cite[Fig.~2]{app}: the empty tile (\texttt{.}),
the crossing (\texttt{+}), vertical and horizontal pipes (\texttt{V},
\texttt{H}), the elbow joining south to east (\texttt{r}), and the elbow
joining west to north (\texttt{j}), such that a pipe enters at every south
edge, one leaves at every east edge, and no pipe meets the north or west
boundary. It is \emph{reduced} if every two pipes cross at most once. The
\(2\times2\) flips of \cite[Fig.~6]{app} replace a \(2\times2\) window by the
other tiling with the same boundary edges, when there is one. Let
\(\Gamma_n\) be the graph whose vertices are the reduced bumpless pipe dreams
(RBPDs) of size \(n\), with an edge for each flip joining two RBPDs. Anderson,
Panova and Petrov \cite[Conjecture~4.2]{app} conjectured:
\begin{quote}
For every \(n\), the set of reduced bumpless pipe dreams of size \(n\) is
connected under the \(2\times 2\) flips (Figure~6) that preserve reducedness.
\end{quote}
Thus \(\Gamma_n\) contains RBPDs of all boundary permutations, and a flip
may change the permutation. The authors verified the conjecture for
\(n\leq 8\) \cite{app}. (Section and figure numbers refer to the arXiv version
of \cite{app}; the published PNAS version states Conjecture~4.2 in the same
words.)

\begin{theorem}\label{thm:main}
\(\Gamma_9\) contains two disjoint sets of \(16\) vertices, each closed under
all Figure~6 flips that preserve reducedness (no such flip leads out of the
set), and neither containing the identity Rothe RBPD. In particular
Conjecture~4.2 is false.
\end{theorem}

\paragraph{The certificate.}
The following grid, with rows listed from top to bottom, is a reduced
bumpless pipe dream with \(16\) crossings, the minimum number in either
set:
\begin{verbatim}
....rHHHH
....VrHHH
....VVrHH
..rH++jrH
..Vr+jr+H
rH+++H++H
Vr+++HjVr
VVVVVrH++
VVVVVVr++
\end{verbatim}
To test reducedness, start every pipe at its south boundary edge and
follow it through the tiles, continuing straight at crossings. At each
crossing record the unordered pair of starting columns; no pair may
repeat. Applying every flip of Figure~6 at each of the \(64\) windows and
keeping the outcomes that pass this test, breadth-first search from the
displayed grid reaches exactly \(16\) RBPDs, joined by \(24\) edges. It
never reaches the identity Rothe RBPD, and all \(170\) flips that leave
this set produce non-reduced pipe dreams. Transposing the grid and
interchanging \texttt{H} and \texttt{V} reverses the paths and preserves
both reducedness and the set of Figure~6 flips; it gives a second set of
\(16\) RBPDs with the same properties, disjoint from the first. This proves
Theorem~\ref{thm:main}; the \(32\) grids are listed in the ancillary files
\texttt{component1.txt} and \texttt{component2.txt}, and
\texttt{verify.py} repeats the checks above using only the Python standard
library.

\paragraph{Permutations.}
We label an RBPD by the permutation \(w\) for which \(w(i)\) is the south
entrance column of the pipe leaving at east row \(i\); this is the convention
implied by the Rothe diagram and Rothe RBPD of \cite[Sect.~4.1 and Fig.~4]{app}.
The displayed grid has \(w=543821967\). The proofs of Lemma~5.3 and
Proposition~6.3 in \cite{app} use the inverse labeling, in which \(w(k)\)
is the exit row of the pipe entering at column \(k\); in that labeling the
displayed grid has \(w=653218947\). Each set contains eight permutations
(Table~\ref{tab:perms}); the permutations of the second set are the
inverses of those of the first, as transposition predicts. No Rothe RBPD of
any of these permutations lies in either set, so each of the sixteen
permutations has many RBPDs outside the certificate.

\begin{table}[htbp]
\centering
\small
\begin{tabular}{r r r@{\qquad}r r r}
\multicolumn{3}{c}{First set} & \multicolumn{3}{c}{Transposed set}\\
\hline
\(w\) & In set & Total & \(w\) & In set & Total\\
\hline
543821967 & 1 & 140 & 653218947 & 1 & 140\\
543821976 & 1 & 525 & 653219847 & 1 & 525\\
548321967 & 3 &  98 & 654218937 & 3 &  98\\
548321976 & 3 & 385 & 654219837 & 3 & 385\\
584321967 & 3 &  56 & 654318927 & 3 &  56\\
584321976 & 3 & 231 & 654319827 & 3 & 231\\
854321967 & 1 &  21 & 654328917 & 1 &  21\\
854321976 & 1 &  91 & 654329817 & 1 &  91\\
\hline
\end{tabular}
\caption{The boundary permutations in each set, the number of RBPDs of
the set having that permutation, and the total number of order-\(9\)
RBPDs with that permutation. The totals equal the principal
specializations \(\mathfrak S_w(1,\dots,1)\); we computed each by the
descent recurrence \cite[Sect.~2.2]{app} and by direct enumeration of all
RBPDs of \(w\) from its Rothe RBPD by moves preserving the permutation
(cf.\ the droops of \cite{lls,app}), and the two counts agree.}
\label{tab:perms}
\end{table}

\section{Minimality and the full census}

We enumerated all RBPDs of size \(n\leq 9\) independently of the flips. Read
the rows from south to north, starting with all \(n\) columns occupied, and
encode the occupied vertical edges below each row by the set of their
column indices. For a candidate next set, the horizontal edge east of
column \(j\) is occupied exactly when the difference between the cumulative
south and north occupancies through \(j\) is one; a row is legal only when
all these differences are zero or one and the last is one. The row is then
determined, and we follow the pipe labels through it and discard any branch
at its first repeated crossing pair. The totals for sizes \(1\) through \(9\)
are
\[
1,\ 2,\ 7,\ 41,\ 393,\ 6{,}080,\ 150{,}371,\ 5{,}903{,}710,\ 365{,}973{,}851,
\]
in agreement with the sum of \(\mathfrak S_w(1,\dots,1)\) over \(w\in S_n\),
that is, with OEIS A331920 \cite{oeis}. On the set of RBPDs, a flip swaps two
adjacent unequal bits of an internal vertical-edge set, provided both
adjacent rows remain legal and the result is reduced. Of the \(6^4\) tile
assignments of a \(2\times2\) window, \(82\) have matching internal edges.
They have \(66\) distinct boundary patterns, of which \(50\) admit a unique
filling and \(16\) admit two; the latter are exactly the \(16\) flip pairs of
Figure~6.

Union--find over this edge set gives the following (two separately written
programs agree).
\begin{proposition}
\(\Gamma_n\) is connected for \(1\leq n\leq 8\). \(\Gamma_9\) has exactly
three components, of sizes \(16\), \(16\) and \(365{,}973{,}819\).
\end{proposition}
Hence \(9\) is the smallest order at which Conjecture~4.2 fails. For
\(n=8\) the program also finds exactly ten non-identity RBPDs with no flip
increasing the height function, matching the ten ``stuck'' RBPDs of
\cite[Sect.~5.1]{app}. The proposition is a computer result; the
counterexample itself rests only on the certificate of
Theorem~\ref{thm:main}, which does not need the census.

\section{Larger orders}

Prepend a blank column to the order-\(9\) grid and append a bottom row
consisting of an \texttt{r} and one \texttt{+} per old column, so that every
old pipe is crossed once by the new pipe. Starting from the \(20\)-crossing
RBPD of the first set (the one with maximal crossing number),
we apply this padding \(1\), \(2\) and \(3\) times to obtain a seed of
order \(10\), \(11\), \(12\). Breadth-first search from the seed, applying
all Figure~6 flips and keeping the reduced outcomes, terminates with a
component of size \(256\), \(5{,}376\) and \(162{,}048\), respectively, none
containing the identity. The order-\(12\) seed (an RBPD with \(50\) crossings)
is in \texttt{seed-n12.txt}, and \texttt{python3 verify.py seed-n12.txt}
repeats the search. These are component censuses from specific seeds, not
censuses of \(\Gamma_n\). The padded copies of the \(16\) RBPDs of order
\(9\) are not themselves a closed set: new flips at the seam leave it. We do
not claim that \(\Gamma_n\) is disconnected for every \(n\geq 9\).

\section{Prior work}

The flips of Figure~6 appear in the full alternating sign matrix lattice in
the work of Escobar, Klein and Weigandt \cite{ekw}, which does not address
connectivity of the reduced graph. Anderson, Panova and Petrov establish
connectivity once the flips are supplemented by droops and undroops
\cite[Sect.~5.2]{app}; that theorem, and the MCMC sampler built on it,
concern a larger move set and are not affected by the counterexample. Up to
28 November 2026 we searched the citing records in Semantic Scholar and
OpenAlex, arXiv, MathOverflow, Mathematics Stack Exchange, GitHub,
Zenodo and general web search, and found no earlier resolution of
Conjecture~4.2; the arXiv record of \cite{app} has only version~1. This
search does not establish priority.

\section*{Acknowledgements}
The authors contributed equally. Computations and parts of the drafting
were assisted by AI tools; the authors verified all results and take
full responsibility for the content.
G.~K.~Chahal thanks her parents, Ranjeet Singh Chahal and Sukhmani Kaur
Chahal, for their constant encouragement. A.~Gupta thanks his parents,
Manish Kumar Gupta and Reena Gupta, and his brother, Mukund Gupta, for
their unwavering support.

\end{document}